\documentclass[aos,preprint]{imsart}

\usepackage{graphicx} 
\usepackage{epsfig}
\usepackage{varioref}
\usepackage{enumitem}
\RequirePackage[OT1]{fontenc}
\RequirePackage{amsmath, latexsym, amsthm, amsfonts, amssymb}
\RequirePackage[authoryear]{natbib}
\RequirePackage[colorlinks,citecolor=blue,urlcolor=blue]{hyperref}

\startlocaldefs
\numberwithin{equation}{section}
\theoremstyle{plain}
\newtheorem{thm}{Theorem}[section]
\newtheorem{prop}{Proposition}[section]

\theoremstyle{definition} 
\newtheorem{exam}{Example}[section]
\newtheorem{subexam}{Subexample}[exam]
\newtheorem{rem}{Remark}[section]

\makeindex

\newcommand{\beao}{\begin{eqnarray*}}
\newcommand{\eeao}{\end{eqnarray*}\noindent}

\newcommand{\beam}{\begin{eqnarray}}
\newcommand{\eeam}{\end{eqnarray}\noindent}

\newcommand{\beqq}{\begin{equation}}
\newcommand{\eeqq}{\end{equation}\noindent}

\newcommand{\bce}{\begin{center}}
\newcommand{\ece}{\end{center}}

\newcommand{\barr}{\begin{array}}
\newcommand{\earr}{\end{array}}
\endlocaldefs

\begin{document}

\begin{frontmatter}

\title{Semiparametrically Efficient Estimation of Constrained Euclidean Parameters}
\runtitle{Constrained Euclidean Parameters}

\begin{aug}
\author{\fnms{Chris A.J.} \snm{Klaassen}\thanksref{t1}\ead[label=e1]{c.a.j.klaassen@uva.nl}}
\and
\author{\fnms{Nanang} \snm{Susyanto}\thanksref{t1}
\ead[label=e2]{n.susyanto@uva.nl}}

\thankstext{t1}{This research has been supported by the Netherlands Organisation for Scientific Research (NWO) via the project Forensic Face Recognition, 727.011.008.}
\runauthor{Klaassen and Susyanto}

\affiliation{University of Amsterdam}
\address{C.A.J. Klaassen\\
N. Susyanto\\
Korteweg-de Vries Institute for Mathematics \\
University of Amsterdam\\
P.O. Box 94248, 1090 GE Amsterdam, The Netherlands\\
\printead{e1}\\
\phantom{E-mail:\ }\printead*{e2}}
\end{aug}

\begin{abstract}
Consider a quite arbitrary (semi)parametric model with a Euclidean parameter of interest and assume that an asymptotically (semi)parametrically efficient estimator of it is given. If the parameter of interest is known to lie on a general surface (image of a continuously differentiable vector valued function), we have a submodel in which this constrained Euclidean parameter may be rewritten in terms of a lower-dimensional Euclidean parameter of interest. An estimator of this underlying parameter is constructed based on the original estimator, and it is shown to be (semi)parametrically efficient. It is proved that the efficient score function for the underlying parameter is determined by the efficient score function for the original parameter and the Jacobian of the function defining the general surface, via a chain rule for score functions. Efficient estimation of the constrained Euclidean parameter itself is considered as well.

Our general estimation method is applied to location-scale, Gaussian copula and semiparametric regression models, and to parametric models under linear restrictions.
\end{abstract}



\end{frontmatter}

\section{Introduction}\label{I}

Let $X_1,\dots,X_n$ be i.i.d. copies of $X$ taking values in the measurable space $({\cal X},{\cal A})$ in a semiparametric model with Euclidean parameter $\theta \in \Theta$ where $\Theta$ is an open subset of $\mathbb{R}^k.$ We denote this semiparametric model by
\begin{equation}\label{model}
{\cal P}=\left\{P_{\theta, G}\ :\ \theta \in \Theta,\  G \in {\cal
G} \right\}.
\end{equation}
Typically, the nuisance parameter space $\cal G$ is a subset of a
Banach or Hilbert space. This space may also be finite dimensional, thus resulting in a parametric model.

We assume an asymptotically efficient estimator $\hat \theta_n=\hat \theta_n(X_1,\dots,X_n)$ is given of the parameter of interest $\theta,$ which under regularity conditions means that
\begin{equation}\label{defefficientestimator}
\sqrt n\left(\hat \theta_n-\theta-\frac1n \sum\limits_{i=1}^n {\tilde \ell}(X_i;\theta,G,{\cal P})\right)\rightarrow_{P_{\theta,G}} 0
\end{equation}
holds. Here ${\tilde \ell}(\cdot;\theta,G,{\cal P})$ is the efficient influence function at $P_{\theta,G}$ for estimation of $\theta$ within $\cal P$ and
\begin{equation}\label{scoretheta}
{\dot \ell}(\cdot; \theta,G,{\cal P}) = \left(\int_{\cal X}{\tilde \ell}(x;\theta,G,{\cal P}){\tilde \ell}^T\!(x;\theta,G,{\cal P})dP_{\theta,G}(x) \right)^{-1} {\tilde \ell}(\cdot;\theta,G,{\cal P})
\end{equation}
is the corresponding efficient score function at $P_{\theta,G}$ for estimation of $\theta$ within $\cal P.$

The topic of this paper is asymptotically efficient estimation when it is known that $\theta$ lies on a general surface, or equivalently, when it is known that $\theta$ is determined by a lower dimensional parameter via a continuously differentiable function, which we denote by
\begin{equation}\label{flat}
\theta = f(\nu), \quad \nu \in N.
\end{equation}
Here $f:N\subset\mathbb{R}^d\to\mathbb{R}^k$ with $d<k$ is known, $N$ is open, the Jacobian
 \begin{equation}\label{Jacobian}
 {\dot f}(\nu) = \left( \frac{\partial f_i(\nu)}{\partial \nu_j} \right)_{j=1,\dots,d}^{i=1,\dots,k}
 \end{equation}
 of $f$ is assumed to be of full rank on $N,$ and $\nu$ is the unknown $d$-dimensional parameter to be estimated.
Thus, we focus on the (semi)parametric model
\begin{equation}\label{submodel}
{\cal Q}=\left\{P_{f(\nu), G}\ :\ \nu \in N,\ G \in {\cal G} \right\} \subset {\cal P}.
\end{equation}
The first main result of this paper is that a semiparametrically efficient estimator of $\nu,$ the parameter of interest, has to be asymptotically linear with efficient score function for estimation of $\nu$ equal to
\begin{equation}\label{scorenu}
{\dot \ell}(\cdot; \nu,G,{\cal Q}) = {\dot f}^T\!(\nu) {\dot \ell}(\cdot; \theta,G,{\cal P}).
\end{equation}
Such a semiparametrically efficient estimator of the parameter of interest can be defined in terms of $f(\cdot)$ and the efficient estimator $\hat \theta_n$ of $\theta;$ see equation (\ref{efficientestimatornu}) in Section \ref{EEPI}. This is our second main result.
How (\ref{scorenu}) is related to the chain rule for differentiation will be explained in Section \ref{CRSF}, which proves this chain rule for score functions.
The semiparametric lower bound for estimators of $\nu$ is obtained via the H\'ajek-LeCam Convolution Theorem for regular parametric models and without projection techniques in Section \ref{CTMR}. In Section \ref{EEPI} efficient estimators within $\cal Q$ of $\nu$ and $\theta$ are constructed, as well as efficient estimators of $\theta$ under linear restrictions on $\theta.$
The generality of our approach facilitates the analysis of numerous statistical models. We discuss some of such parametric and semiparametric models and related literature in Section \ref{E}. One of the proofs will be given in Appendix \ref{TP}.

The topic of this paper should not be confused with estimation of the parameter $\theta$ when it is known to lie in a subset of the original parameter space described by linear inequalities. A comprehensive treatment of such estimation problems may be found in \cite{vanEeden06}. Our model $\cal Q$ with its constrained Euclidean parameters also differs from the constraint defined models as studied by Bickel et al. (1993, 1998) (henceforth called BKRW), which are defined by restrictions on the distributions in $\cal P.$

\section{The Chain Rule for Score Functions}\label{CRSF}

The basic building block for the asymptotic theory of semiparametric models as presented in e.g. \cite{BKRW93} is the concept of regular parametric model. Let ${\cal P}_\Theta = \left\{ P_\theta\,:\, \theta \in \Theta \right\}$ with $\Theta \subset {\mathbb R}^k$ open be a parametric model with all $P_\theta$ dominated by a $\sigma$-finite measure $\mu$ on $\left({\cal X},{\cal A}\right).$ Denote the density of $P_\theta$ with respect to $\mu$ by $p(\theta)=p(\cdot; \theta, {\cal P}_\Theta)$ and the $L_2(\mu)$-norm by $\parallel \cdot \parallel_\mu.$ If for each $\theta_0 \in \Theta$ there exists a $k$-dimensional column vector ${\dot \ell}(\theta_0, {\cal P}_\Theta)$ of elements of $L_2(P_{\theta_0}),$ the so-called score function, such that the Fr\'echet differentiability
\begin{eqnarray}\label{Frechet}
\lefteqn{\parallel \sqrt{p(\theta)} - \sqrt{p(\theta_0)} -\tfrac 12 \left(\theta-\theta_0 \right)^T{\dot \ell}(\theta_0,{\cal P}_\Theta) \sqrt{p(\theta_0)} \parallel_\mu \nonumber} \\
&& \qquad  \qquad \qquad \qquad \qquad = o(|\theta -\theta_0|),\quad \theta \to \theta_0,
\end{eqnarray}
holds and the $k \times k$ Fisher information matrix
\begin{equation}\label{Fisher}
I(\theta_0) = \int_{\cal X} {\dot \ell}(\theta_0,{\cal P}_\Theta) {\dot \ell}^T\!(\theta_0,{\cal P}_\Theta) dP_{\theta_0}
\end{equation}
is nonsingular, and, moreover, the map $\theta \mapsto {\dot \ell}(\theta,{\cal P}_\Theta) \sqrt{p(\theta)}$ from $\Theta$ to $L_2^k(\mu)$ is continuous, then ${\cal P}_\Theta$ is called a {\em regular parametric} model. Often the score function may be determined by computing the logarithmic derivative of the density with respect to $\theta;$ cf. Proposition 2.1.1 of \cite{BKRW93}.
We will call $\cal P$ from (\ref{model}) a {\em regular semiparametric} model if for all $G \in {\cal G}$
\begin{equation}\label{submodel1}
{\cal P}_{\Theta,G} = \left\{P_{\theta, G}\,:\, \theta \in \Theta \right\}
\end{equation}
is a regular parametric model.

Fix $\theta_0 \in \Theta$ and $G_0 \in {\cal G},$ and write $P_{\theta_0, G_0}=P_0.$
Let $\psi\,:\,\Theta \to {\cal G}$ with $\psi(\theta_0)=G_0$ be such that
\begin{equation}\label{submodelpsi}
{\cal P}_\psi = \left\{ P_{\theta, \psi(\theta)} \,:\, \theta \in \Theta \right\}
\end{equation}
is a regular parametric submodel of $\cal P$ with score function ${\dot \ell}(\theta_0, {\cal P}_\psi)$ at $\theta_0$ and Fisher information matrix $I(\theta_0,{\cal P}_\psi),$ say.
Let the density of $P_{\theta, \psi(\theta)}$ with respect to $\mu$ be denoted by $q(\theta).$ Since ${\cal P}_\psi$ is a regular parametric model the score function ${\dot \ell}(\theta_0, {\cal P}_\psi)$ for $\theta$ at $\theta_0$ within ${\cal P}_\psi$ satisfies (cf. (\ref{Frechet}))
\begin{eqnarray}\label{Frechetefficient}
\lefteqn{\parallel \sqrt{q(\theta)} - \sqrt{q(\theta_0)} -\tfrac 12 \left(\theta -\theta_0 \right)^T{\dot \ell}(\theta_0,{\cal P}_\psi)\sqrt{q(\theta_0)} \parallel_\mu \nonumber } \\
&& \qquad \qquad \qquad \qquad \qquad \qquad = o(|\theta-\theta_0|),\quad \theta \to \theta_0.
\end{eqnarray}

Considering now the (semi)parametric submodel $\cal Q$ from(\ref{submodel}) we fix $\nu_0$ and write $f(\nu_0)=\theta_0$ and $f(\nu) = \theta.$ Within $\cal Q$ the Fr\'echet differentiability (\ref{Frechetefficient}) yields
\begin{eqnarray}\label{FrechetefficientQ}
\lefteqn{\parallel \sqrt{q(f(\nu))} - \sqrt{q(f(\nu_0))}
-\tfrac 12 \left(f(\nu)- f(\nu_0)\right)^T{\dot \ell}(f(\nu_0),{\cal P}_\psi) \sqrt{q(f(\nu_0))} \parallel_\mu \nonumber}\\
&& \qquad \qquad \qquad \qquad \qquad \qquad \qquad = o(|f(\nu)-f(\nu_0)|), \quad f(\nu) \to f(\nu_0),
\end{eqnarray}
and hence
\begin{eqnarray}\label{FrechetefficientQ2}
\lefteqn{\parallel \sqrt{q(f(\nu))} - \sqrt{q(f(\nu_0))}
-\tfrac 12 (\nu - \nu_0)^T {\dot f}^T\!(\nu_0) {\dot \ell}(\theta_0,{\cal P}_\psi)\sqrt{q(f(\nu_0))} \parallel_\mu \nonumber}\\
&& \qquad \qquad \qquad \qquad \qquad \qquad \qquad \qquad \qquad = o(|\nu - \nu_0|), \quad \nu \to \nu_0,
\end{eqnarray}
in view of the differentiability of $f(\cdot).$ Since ${\dot f}(\cdot)$ is continuous, this means that
\begin{equation}\label{Qregularsubmodel}
{\cal Q}_\psi = \left\{ P_{f(\nu),\psi(f(\nu))}\,:\, \nu \in N \right\}
\end{equation}
is a regular parametric submodel of $\cal Q$ with score function
\begin{equation}\label{chainscore}
{\dot \ell}(\nu_0,{\cal Q}_\psi) = {\dot f}^T\!(\nu_0){\dot \ell}(\theta_0,{\cal P}_\psi)
\end{equation}
for $\nu$ at $P_0$ and Fisher information matrix
\begin{equation}\label{FisherQ}
{\dot f}^T\!(\nu_0) I(\theta_0,{\cal P}_\psi){\dot f}(\nu_0)
= {\dot f}^T\!(\nu_0)\int_{\cal X} {\dot \ell}(\theta_0,{\cal P}_\psi) {\dot \ell}^T\!(\theta_0,{\cal P}_\psi)dP_0\ {\dot f}(\nu_0) .
\end{equation}
We have proved
\begin{prop}\label{chaindifferentiability}
Let $\cal P$ as in (\ref{model}) be a regular semiparametric model and let $\cal Q$ as in (\ref{submodel}) be a regular semiparametric submodel with $f(\cdot)$ and ${\dot f}(\cdot)$ defined as in and below (\ref{flat}) and (\ref{Jacobian}). If there exists a regular parametric submodel ${\cal P}_\psi$ of $\cal P$ with score function ${\dot \ell}(\theta_0, {\cal P}_\psi)$ for $\theta$ at $\theta_0=f(\nu_0),$ then there exists a regular parametric submodel ${\cal Q}_\psi$ of $\cal Q$ with score function ${\dot \ell}(\nu_0,{\cal Q}_\psi)$ for $\nu$ at $\nu_0$ satisfying (\ref{chainscore}).
\end{prop}
This Proposition is also valid for parametric models, as may be seen by choosing $\cal G$ finite dimensional or even degenerate. The basic version of the chain rule for score functions is for such a parametric model ${\cal P}_\Theta.$ We have chosen the more elaborate formulation of Proposition \ref{chaindifferentiability} since we are going to apply the chain rule for such parametric submodels ${\cal P}_\psi$ of semiparametric models $\cal P.$

\section{Convolution Theorem and Main Result}\label{CTMR}

An estimator ${\hat \theta}_n$ of $\theta$ within the regular semiparametric model $\cal P$ is called (locally) regular at $P_0=P_{\theta_0,G_0}$ if it is (locally) regular at $P_0$ within ${\cal P}_\psi$ for all regular parametric submodels ${\cal P}_\psi$ of $\cal P$ containing $P_{\Theta,G_0}.$
According to the H\'ajek-LeCam Convolution Theorem for regular parametric models (see e.g. Section 2.3 of \cite{BKRW93}) this implies that such a regular estimator ${\hat \theta}_n$ of $\theta$ within $\cal P$ has a limit distribution under $P_0$ that is the convolution of a normal distribution with mean 0 and covariance matrix $I^{-1}(\theta_0,{\cal P}_\psi)$ and another distribution, for any regular parametric submodel ${\cal P}_\psi$ containing $P_0.$ If there exists $\psi=\psi_0$ such that this last distribution is degenerate at 0, we call ${\hat \theta}_n$ (locally) efficient at $P_0$ and ${\cal P}_{\psi_0}$ a least favorable parametric submodel for estimation of $\theta$ within $\cal P$ at $P_0.$
Then the H\'ajek-LeCam Convolution Theorem also implies that ${\hat \theta}_n$ is asymptotically linear in the efficient influence function ${\tilde \ell}(\theta_0,G_0, {\cal P}) = {\tilde \ell}(\cdot; \theta_0,G_0, {\cal P})$ satisfying
\begin{equation}\label{efficientinfluencefunction}
{\tilde \ell}(\theta_0, G_0, {\cal P}) = {\tilde \ell}(\theta_0, {\cal P}_{\psi_0}) = I^{-1}(\theta_0,{\cal P}_{\psi_0}){\dot \ell}(\theta_0, {\cal P}_{\psi_0}),
\end{equation}
which means
\begin{equation}\label{efficientestimatorunderpsi0}
\sqrt n\left({\hat \theta}_n-\theta_0-\frac1n \sum\limits_{i=1}^n {\tilde \ell}(X_i; \theta_0,G_0,{\cal P})\right)\rightarrow_{P_0} 0.
\end{equation}

The argument above can be extended to the more general situation that there exists a least favorable sequence of parametric submodels indexed by $\psi_j\,,j=1, 2, \dots,$ such that the corresponding score functions ${\dot \ell}(\theta_0,{\cal P}_{\psi_j})$ for $\theta$ at $\theta_0$ within model ${\cal P}_{\psi_j}$ converge in $L_2^k(P_0)$ to ${\dot \ell}(\theta_0,G_0,{\cal P}) = {\dot \ell}(\cdot;\theta_0,G_0,{\cal P}),$ say. A regular estimator ${\hat \theta}_n$ of $\theta$ within $\cal P$ is called efficient then, if it is asymptotically linear as in (\ref{efficientestimatorunderpsi0}) with efficient influence function ${\tilde \ell}(\theta_0,G_0, {\cal P}) = {\tilde \ell}(\cdot;\theta_0,G_0, {\cal P})$ satisfying
\begin{eqnarray}\label{efficientinfluencefunctiongeneral}
\lefteqn{{\tilde \ell}(\theta_0,G_0, {\cal P})= \left(\int_{\cal X}{\dot \ell}(\theta_0,G_0, {\cal P}){\dot \ell}^T(\theta_0,G_0, {\cal P}) dP_0 \right)^{-1}
{\dot \ell}(\theta_0,G_0, {\cal P})\nonumber} \\
&& \qquad \qquad \qquad \qquad \qquad \qquad = I^{-1}(\theta_0,G_0, {\cal P}){\dot \ell}(\theta_0,G_0, {\cal P}).
\end{eqnarray}

Indeed, by the Convolution Theorem for regular parametric models the convergence
\begin{equation}\label{convolution}
\begin{pmatrix}
\sqrt n\left({\hat \theta}_n-\theta_0-\frac1n \sum\limits_{i=1}^n {\tilde \ell}(X_i; \theta_0,{\cal P}_{\psi_j})\right)\\
\frac1{\sqrt n} \sum\limits_{i=1}^n {\tilde \ell}(X_i; \theta_0,{\cal P}_{\psi_j})
\end{pmatrix}
\rightarrow_{P_0}
\begin{pmatrix}
R_j \\ Z_j
\end{pmatrix}
\end{equation}
holds with the $k$-vectors $R_j$ and $Z_j$ independent and $Z_j$ normal with mean 0 and covariance matrix $I^{-1}(\theta_0,{\cal P}_{\psi_j}).$ Taking limits as $j \to \infty$ we see by tightness arguments and by the convergence of ${\dot \ell}(\theta_0,{\cal P}_{\psi_j})$ to ${\dot \ell}(\theta_0,G_0, {\cal P})$ in $L_2^k(P_0),$ that also
\begin{equation}\label{convolution2}
\begin{pmatrix}
\sqrt n\left({\hat \theta}_n-\theta_0-\frac1n \sum\limits_{i=1}^n {\tilde \ell}(X_i; \theta_0,G_0, {\cal P})\right)\\
\frac1{\sqrt n} \sum\limits_{i=1}^n {\tilde \ell}(X_i; \theta_0,G_0, {\cal P})
\end{pmatrix}
\rightarrow_{P_0}
\begin{pmatrix}
R_{\cal P} \\ Z_{\cal P}
\end{pmatrix}
\end{equation}
holds with $R_{\cal P}$ and $Z_{\cal P}$ independent. If $R_{\cal P}$ is degenerate at 0, then ${\hat \theta}_n$ is locally asymptotically efficient at $P_0$ within $\cal P$ and the sequence of regular parametric submodels ${\cal P}_{\psi_j}$ is least favorable indeed.

Now, let us assume such a least favorable sequence and efficient estimator ${\hat \theta}_n$ exist at $P_0=P_{\theta_0,G_0}$ with $\theta_0=f(\nu_0)$ and $f(\cdot)$ from (\ref{flat}) and (\ref{Jacobian}) continuously differentiable. By the chain rule for score functions from Proposition \ref{chaindifferentiability} the score function ${\dot \ell}(\nu_0,{\cal Q}_{\psi_j})$ for $\nu$ at $\nu_0$ within ${\cal Q}_{\psi_j}$ satisfies
\begin{equation}\label{chainscore2}
{\dot \ell}(\nu_0,{\cal Q}_{\psi_j}) = {\dot f}^T\!(\nu_0){\dot \ell}(\theta_0,{\cal P}_{\psi_j})
\end{equation}
and hence the corresponding influence function ${\tilde \ell}(\nu_0,{\cal Q}_{\psi_j})$ satisfies
\begin{equation}\label{influencefunctionQpsij}
{\tilde \ell}(\nu_0,{\cal Q}_{\psi_j})=
\left({\dot f}^T\!(\nu_0) I(\theta_0,{\cal P}_{\psi_j}){\dot f}(\nu_0)\right)^{-1}{\dot f}^T\!(\nu_0){\dot \ell}(\theta_0,{\cal P}_{\psi_j}).
\end{equation}
Let ${\hat \nu}_n$ be a locally regular estimator of $\nu$ at $P_0$ within the regular semiparametric model $\cal Q.$
By the convergence of ${\dot \ell}(\theta_0,{\cal P}_{\psi_j})$ to ${\dot \ell}(\theta_0,G_0, {\cal P})$ in $L_2^k(P_0),$ the influence functions from (\ref{influencefunctionQpsij}) converge in $L_2^d(P_0)$ to
\begin{equation}\label{influencefunctionQ}
{\tilde \ell}(\nu_0,G_0,{\cal Q})=
\left({\dot f}^T\!(\nu_0) I(\theta_0,G_0, {\cal P}){\dot f}(\nu_0)\right)^{-1}{\dot f}^T\!(\nu_0){\dot \ell}(\theta_0,G_0, {\cal P})
\end{equation}
and the argument leading to (\ref{convolution2}) yields the convergence
\begin{equation}\label{convolutionnu}
\begin{pmatrix}
\sqrt n\left({\hat \nu}_n-\nu_0-\frac1n \sum\limits_{i=1}^n {\tilde \ell}(X_i; \nu_0,G_0,{\cal Q})\right)\\
\frac1{\sqrt n} \sum\limits_{i=1}^n {\tilde \ell}(X_i; \nu_0,G_0,{\cal Q})
\end{pmatrix}
\rightarrow_{P_0}
\begin{pmatrix}
R_{\cal Q} \\ Z_{\cal Q}
\end{pmatrix}
\end{equation}
with $R_{\cal Q}$ and $Z_{\cal Q}$ independent. Note that $Z_{\cal Q}$ has a normal distribution with mean 0 and covariance matrix
\begin{equation}\label{informationboundnu}
I^{-1}(\nu_0, G_0, {\cal Q}) = \left({\dot f}^T\!(\nu_0) I(\theta_0,G_0, {\cal P}){\dot f}(\nu_0)\right)^{-1}.
\end{equation}
Under an additional condition on $f(\cdot)$ we shall construct an estimator ${\hat \nu}_n$ of $\nu$ based on ${\hat \theta}_n$ for which $R_{\cal Q}$ is degenerate. This construction of ${\hat \nu}_n$ will be given in the next section together with a proof of its efficiency, and this will complete the proof of our main result formulated as follows.

\begin{thm}\label{mainresult}
Let $\cal P$ from (\ref{model}) be a regular semiparametric model with $P_0=P_{\theta_0,G_0} \in {\cal P}, \theta_0=f(\nu_0),$ and $f(\cdot)$ from (\ref{flat}) and (\ref{Jacobian}) continuously differentiable. Furthermore, let $f(\cdot)$ have an inverse on $f(N)$ that is differentiable with a bounded Jacobian. If there exists a least favorable sequence of regular parametric submodels ${\cal P}_{\psi_j}$ and an asymptotically efficient estimator ${\hat \theta}_n$ of $\theta$ satisfying
(\ref{convolution2}) with $R_{\cal P}=0$ a.s., then there exists a least favorable sequence of regular parametric submodels ${\cal Q}_{\psi_j}$ of the restricted model $\cal Q$ from (\ref{submodel}) and an asymptotically efficient estimator ${\hat \nu}_n$ of $\nu$ satisfying (\ref{convolutionnu}) with $R_{\cal Q}=0$ a.s. and attaining the asymptotic information bound (\ref{informationboundnu}).
\end{thm}
Note that the convolution result (\ref{convolutionnu}) and (\ref{influencefunctionQ}) also holds if the convergent sequence of regular parametric submodels ${\cal P}_{\psi_j}$ is not least favorable, and that it implies by the central limit theorem that the limit distribution of ${\sqrt n}\left({\hat \nu}_n -\nu_0 \right)$ is the convolution of a normal distribution with mean 0 and covariance matrix
\begin{equation}\label{boundcovariancematrix}
I^{-1}(\nu_0, G_0, {\cal Q}) = \left({\dot f}^T\!(\nu_0) I(\theta_0,G_0, {\cal P}){\dot f}(\nu_0)\right)^{-1}
\end{equation}
and the distribution of $R_{\cal Q}.$

\section{Efficient Estimator of the Parameter of Interest}\label{EEPI}
There are many ways of constructing efficient estimators in (semi)parametric models. One of the common approaches is upgrading a $\sqrt{n}$-consistent estimator as in Sections 2.5 and 7.8 of \cite{BKRW93}. A somewhat different upgrading approach is used in the following construction.
\begin{thm}\label{onestep}
Consider the situation of Theorem \ref{mainresult}. If the symmetric positive definite $k \times k$-matrix ${\hat I}_n$ is a consistent estimator of $I(\theta, G, {\cal P})$ within $\cal P$ and ${\bar \nu}_n$ is a $\sqrt n$-consistent estimator of $\nu$ within $\cal Q,$ then
\begin{equation}\label{efficientestimatornu}
{\hat \nu}_n = {\bar \nu}_n + \left({\dot f}^T\!({\bar \nu}_n) {\hat I}_n{\dot f}({\bar \nu}_n)\right)^{-1}
{\dot f}^T\!({\bar \nu}_n) {\hat I}_n \left[{\hat \theta}_n - f\left({\bar \nu}_n \right) \right]
\end{equation}
is efficient, i.e., it satisfies (\ref{convolutionnu}) with $R_{\cal Q}=0$ a.s.
\end{thm}
\noindent {\bf Proof}
The continuity of ${\dot f}(\cdot)$ and the consistency of ${\bar \nu}_n$ and ${\hat I}_n$ imply that
\begin{equation}\label{Kn}
{\hat K}_n = \left({\dot f}^T\!({\bar \nu}_n) {\hat I}_n{\dot f}({\bar \nu}_n)\right)^{-1} {\dot f}^T\!({\bar \nu}_n) {\hat I}_n
\end{equation}
converges in probability under $P_0$ to
\begin{equation}\label{K0}
K_0 = \left({\dot f}^T\!(\nu_0) I(\theta_0,G_0,{\cal P}){\dot f}(\nu_0)\right)^{-1} {\dot f}^T\!(\nu_0) I(\theta_0,G_0,{\cal P}).
\end{equation}
This means that ${\hat K}_n$ consistently estimates $K_0.$ In view of (\ref{efficientestimatornu}), (\ref{influencefunctionQ}), (\ref{efficientinfluencefunctiongeneral}), and (\ref{convolution2}) with $R_{\cal P}=0$ we obtain
\begin{eqnarray}\label{efficiencyargument}
\lefteqn{{\sqrt n}\left({\hat \nu}_n - \nu_0 -\frac1n \sum\limits_{i=1}^n {\tilde \ell}(X_i; \nu_0, G_0,{\cal Q})\right)\nonumber }\\
&& = {\sqrt n}\left({\bar \nu}_n - \nu_0 + {\hat K}_n\left[{\hat \theta}_n - f\left({\bar \nu}_n \right) \right]
-\frac1n \sum\limits_{i=1}^n K_0{\tilde \ell}(X_i;\theta_0,G_0,{\cal P}) \right) \nonumber \\
&& = {\sqrt n}\left({\bar \nu}_n - \nu_0 - {\hat K}_n\left[f\left({\bar \nu}_n \right) - f(\nu_0) \right] \right) \nonumber \\
&& \qquad \qquad \qquad \qquad + \left[{\hat K}_n - K_0\right] \frac1{\sqrt n} \sum\limits_{i=1}^n {\tilde \ell}(X_i;\theta_0,G_0,{\cal P}) + o_p(1).
\end{eqnarray}
By the consistency of ${\hat K}_n$ the second term at the right hand side of (\ref{efficiencyargument}) converges to 0 in probability under $P_0$ in view of the central limit theorem. Because $f\left({\bar \nu}_n \right) = f(\nu_0) + {\dot f}(\nu_0)\left({\bar \nu}_n - \nu_0 \right) + o_p\left({\bar \nu}_n - \nu_0 \right)$ holds and $K_0{\dot f}(\nu_0)$ equals the $d \times d$ identity matrix, the first part of the right hand side of (\ref{efficiencyargument}) also converges to 0 in probability under $P_0.$
\hfill$\Box$\\

To complete the proof of Theorem \ref{mainresult} with the help of Theorem \ref{onestep} we will construct a $\sqrt{n}$-consistent estimator ${\bar \nu}_n$ of $\nu$ and subsequently a consistent estimator ${\hat I}_n$ of $I(\theta, G, {\cal P}).$ Let $\parallel \cdot \parallel$ be a Euclidean norm on ${\mathbb R}^k.$ We choose ${\bar \nu}_n$ in such a way that
\begin{equation}\label{triangle}
\parallel f\left({\bar \nu}_n\right) - {\hat \theta}_n \parallel \leq \inf_{\nu \in N} \parallel f(\nu) - {\hat \theta}_n \parallel + \frac 1n
\end{equation}
holds. Of course, if the infimum is attained, we choose ${\bar \nu}_n$ as the minimizer. By the triangle inequality and the $\sqrt n$-consistency of ${\hat \theta}_n$ we obtain
\begin{eqnarray}\label{triangle2}
\lefteqn{\parallel f({\bar \nu}_n) - f(\nu_0) \parallel \leq \inf_{\nu \in N} \parallel f(\nu) - {\hat \theta}_n \parallel + \frac 1n + \parallel f(\nu_0) - {\hat \theta}_n \parallel \nonumber}\\
&& \qquad \qquad \qquad \leq 2\parallel {\hat \theta}_n - f(\nu_0) \parallel +\frac 1n = O_p\left(\frac 1{\sqrt n} \right).
\end{eqnarray}
The assumption from Theorem \ref{mainresult} that $f(\cdot)$ has an inverse on $f(N)$ that is differentiable with a bounded Jacobian, suffices to conclude that (\ref{triangle2}) guarantees $\sqrt n$-consistency of ${\bar \nu}_n.$

In constructing a consistent estimator of the Fisher information matrix based on the given efficient estimator ${\hat \theta}_n,$ we split the sample in blocks as follows. Let $(k_n), (\ell_n),$ and $(m_n)$ be sequences of integers such that $k_n = \ell_n m_n, k_n/n \to \kappa, 0<\kappa<1,$ and $\ell_n \to \infty, m_n \to \infty$ hold as $n \to \infty.$ For $j= 1,\dots,\ell_n$ let ${\hat \theta}_{n,j}$ be the efficient estimator of $\theta$ based on the observations $X_{(j-1)m_n +1}, \dots, X_{jm_n}$ and ${\hat \theta}_{n,0}$ be the efficient estimator of $\theta$ based on the remaining observations $X_{k_n +1}, \dots, X_n.$ Consider the "empirical" characteristic function
\begin{equation}\label{empiricalcharacteristicfunction}
{\hat \phi}_n(t) = \frac 1{\ell_n} \sum _{j=1}^{\ell_n} \exp\left\{ it \sqrt{m_n}\left({\hat \theta}_{n,j} - {\hat \theta}_{n,0} \right) \right\},\ t \in {\mathbb R}^k,
\end{equation}
which we rewrite as
\begin{eqnarray}\label{echfrewritten}
&&{\hat \phi}_n(t) = \exp\left\{ -it \sqrt{m_n}\left({\hat \theta}_{n,0} - \theta_0 \right) \right\} \frac 1{\ell_n} \sum _{j=1}^{\ell_n} \exp\left\{ it \sqrt{m_n}\left({\hat \theta}_{n,j} - \theta_0 \right) \right\} \nonumber \\
&& \qquad \qquad = \exp\left\{ -it \sqrt{m_n}\left({\hat \theta}_{n,0} - \theta_0 \right) \right\} {\tilde \phi}_n(t).
\end{eqnarray}
In view of $m_n/(n - k_n) \to 0$ and (\ref{convolution2}) with $R_{\cal P}=0$ a.s. we see that the first factor at the right hand side of (\ref{echfrewritten}) converges to 1 as $n\to \infty.$ The efficiency of ${\hat \theta}_n$ in (\ref{convolution2}) with $R_{\cal P}=0$ a.s. also implies
\begin{eqnarray}\label{mean}
&& E\left({\tilde \phi}_n(t) \right) = E\left( \exp\left\{ it \sqrt{m_n}\left({\hat \theta}_{n,1} - \theta_0 \right) \right\} \right) \nonumber \\
&& \qquad \qquad \qquad \to E\left( \exp\left\{it Z_{\cal P}\right\} \right)
\end{eqnarray}
as $n\to \infty,$ with $Z_{\cal P}$ normally distributed with mean 0 and covariance matrix $I^{-1}(\theta_0, G_0,{\cal P}).$ Some computation shows
\begin{eqnarray}\label{variancephitilde}
&& E\left(\left| {\tilde \phi}_n(t) - E\left({\tilde \phi}_n(t)\right) \right|^2 \right) \nonumber \\
 && \qquad = \frac 1{\ell_n} \left(1 - \left| E\left(\exp\left\{ it \sqrt{m_n}\left({\hat \theta}_{n,1} - \theta_0 \right) \right\}\right) \right|^2 \right) \leq \frac 1{\ell_n}.
\end{eqnarray}
It follows by Chebyshev's inequality that ${\tilde \phi}_n(t)$  and hence ${\hat \phi}_n(t)$ converges under $P_0 =P_{\theta_0,G_0}$ to the characteristic function of $Z_{\cal P}$ at $t,$
\begin{equation}\label{convergencephihat}
{\hat \phi}_n(t) \rightarrow_{P_0}   E\left( \exp\left\{it Z_{\cal P}\right\} \right) =\exp\left\{- \tfrac 12 t^T I^{-1}(\theta_0, G_0, {\cal P}) t \right\}.
\end{equation}
For every $t \in {\mathbb R}^k$ we obtain
\begin{equation}\label{convergencelogphi}
-2 \log\left(\Re \left({\hat \phi}_n(t)\right)\right) \rightarrow_{P_0} t^T I^{-1}(\theta_0, G_0, {\cal P}) t.
\end{equation}
Choosing $k(k+1)/2$ appropriate values of $t$ we may obtain from (\ref{convergencelogphi}) an estimator of $I^{-1}(\theta_0, G_0, {\cal P})$ and hence of $I(\theta_0, G_0, {\cal P}).$ Indeed, with $t$ equal to the unit vectors $u_i$ we obtain estimators of the diagonal elements of $I^{-1}(\theta_0, G_0, {\cal P})$ and an estimator of its $(i,j)$ element is obtained via $$\log\left(\Re \left({\hat \phi}_n(u_i)\right)\right) + \log\left(\Re \left({\hat \phi}_n(u_j)\right)\right)- \log\left(\Re \left({\hat \phi}_n(u_i + u_j)\right)\right).$$ When needed, the resulting estimator of $I(\theta_0, G_0, {\cal P})$  can be made positive definite by changing appropriate components of it by an asymptotically negligible amount, while the symmetry is maintained.

Under a mild uniform integrability condition it has been shown by \cite{Klaassen87}, that existence of an efficient estimator ${\hat \theta}_n$ of $\theta$ in $\cal P$ implies the existence of a consistent and $\sqrt n$-unbiased estimator of the efficient influence function ${\tilde \ell}(\cdot;\theta, G, {\cal P}).$ Basing this estimator on one half of the sample and taking the average of this estimated efficient influence function at the observations from the other half of the sample, we could have constructed another estimator of the efficient Fisher information. However, this estimator would have been more involved, and, moreover, it needs this extra uniformity condition.

With the help of Theorem \ref{onestep}, the estimator ${\bar \nu}_n$ of $\nu$ from (\ref{triangle}), and the construction via (\ref{convergencelogphi}) of an estimator ${\hat I}_n$ of the efficient Fisher information we have completed our construction of an efficient estimator ${\hat \nu}_n$ as in (\ref{efficientestimatornu}) of $\nu.$ This estimator can be turned into an efficient estimator of $\theta = f(\nu)$ within the model $\cal Q$ from (\ref{submodel}) by
\begin{equation}\label{efficientestimatorfnu}
\tilde\theta_n=f({\hat\nu}_n)
\end{equation}
with efficient influence function
\begin{eqnarray}\label{influencefunctiontheta}
\lefteqn{ {\tilde \ell}(\theta_0,G_0,{\cal Q}) = {\dot f}(\nu_0){\tilde \ell}(\nu_0,G_0,{\cal Q}) \nonumber }\\
&& \qquad = {\dot f}(\nu_0)
\left({\dot f}^T\!(\nu_0) I(\theta_0,G_0, {\cal P}){\dot f}(\nu_0)\right)^{-1}{\dot f}^T\!(\nu_0){\dot \ell}(\theta_0,G_0, {\cal P})
\end{eqnarray}
and asymptotic information bound
\begin{equation}\label{informationboundtheta}
I^{-1}(\theta_0, G_0, {\cal Q}) = {\dot f}(\nu_0)\left({\dot f}^T\!(\nu_0) I(\theta_0,G_0, {\cal P}){\dot f}(\nu_0)\right)^{-1}{\dot f}^T\!(\nu_0).
\end{equation}
Indeed, according to \cite{BKRW93} Section 2.3, $\tilde\theta_n$ is efficient for estimation of $\theta$ under the additional information $\theta=f(\nu).$

\begin{rem}\label{remarklinear}
If $f(\cdot)$ is a linear function, i.e., $\theta = L \nu +\alpha$ holds with the $k \times d$-matrix $L$ of maximum rank $d$, then
\begin{equation}\label{nubarforlinear}
{\bar \nu}_n=(L^TL)^{-1}L^T({\hat \theta}_n -\alpha)
\end{equation}
attains the infimum at the right hand side of (\ref{triangle}).
So, the estimator \eqref{efficientestimatornu} becomes
\begin{equation}\label{efficientnuforlinear}
{\hat \nu}_n = \left(L^T {\hat I}_nL\right)^{-1}
L^T {\hat I}_n\left[\hat \theta_n-\alpha \right]
\end{equation}
with efficient influence function (\ref{influencefunctionQ}) and asymptotic information bound (\ref{informationboundnu}) with ${\dot f}(\nu_0) = L,$
and the estimator from (\ref{efficientestimatorfnu})
\begin{equation}\label{efficientestimatorfnuforlinear}
\tilde\theta_n = L \left(L^T {\hat I}_nL\right)^{-1}
L^T {\hat I}_n\left[\hat \theta_n-\alpha \right] + \alpha.
\end{equation}
Note that ${\tilde \theta}_n$ is the projection of ${\hat \theta}_n$ on the flat $\{\theta \in {\mathbb R}^k\, :\, \theta =L\nu + \alpha, \nu \in {\mathbb R}^d \}$ under the inner product determined by ${\hat I}_n$ (cf. Appendix \ref{TP}) and that the covariance matrix of its limit distribution equals the asymptotic information bound
\begin{equation}\label{informationboundlinear}
I^{-1}(\theta_0, G_0, {\cal Q}) = L \left( L^T I(\theta_0, G_0, {\cal P})L \right)^{-1} L^T.
\end{equation}
Another way to describe this submodel $\cal Q$ with $\theta=L\nu + \alpha$ is by linear restrictions
\begin{equation}\label{restrictions}
{\cal Q} = \left\{P_{L\nu +\alpha}\,:\, \nu \in N, G \in {\cal G} \right\} = \left\{ P_{\theta,G}\,:\, R^T \theta = \beta, \theta \in \Theta, G \in {\cal G} \right\},
\end{equation}
where $R^T\alpha = \beta$ holds and the $k \times d$-matrix $L$ and the $k \times (k-d)$-matrix $R$ are matching such that the columns of $L$ are orthogonal to those of $R$ and the $k \times k$-matrix $(L\, R)$ is of rank $k.$ Note that the open subset $N$ of ${\mathbb R}^d$ determines the open subset $\Theta$ of ${\mathbb R}^k$ and vice versa. See \cite{Cobb28}, \cite{Stone54}, \cite{Nyquist91}, and \cite{Kim95} for some examples of estimation under linear restrictions.

In terms of the restrictions described by $R$ and $\beta$ the efficient estimator ${\tilde \theta}_n$ of $\theta$ from(\ref{efficientestimatorfnuforlinear}) within the submodel $\cal Q$ can be rewritten as
\begin{equation}\label{efficientestimatorfnurestrictions}
{\tilde \theta}_n = {\hat \theta}_n - {\hat I}_n^{-1}R\left(R^T {\hat I}_n^{-1} R\right)^{-1} \left(R^T{\hat \theta}_n - \beta \right),
\end{equation}
with asymptotic information bound
\begin{equation}\label{informationboundrestriction}
L(L^TIL)^{-1}L^T = I^{-1} -I^{-1}R(R^TI^{-1}R)^{-1}R^TI^{-1},\ I = I(\theta_0, G_0, {\cal P}),
\end{equation}
as will be proved in Appendix \ref{TP}.
\end{rem}

\section{Examples}\label{E}
In this section we present five examples, which illustrate our construction of (semi)parametrically efficient estimators. We shall discuss location-scale, Gaussian copula, and semiparametric regression models, and parametric models under linear restrictions.

\begin{exam}\label{locationscale}\textbf{Coefficient of variation known}

Let $g(\cdot)$ be an absolutely continuous density on $({\mathbb R}, \cal{B})$ with mean 0, variance 1, and derivative $g'(\cdot),$ such that $\int [1+x^2] (g'/g(x))^2 g(x) dx$ is finite. Consider the location-scale family corresponding to $g(\cdot).$ Let there be given efficient estimators ${\bar \mu}_n$ and ${\bar \sigma}_n$ of $\mu$ and $\sigma,$ respectively, based on $X_1, \dots, X_n,$ which are i.i.d. with density $\sigma^{-1}g((\cdot - \mu)/\sigma).$ By $I_{ij}$ we denote the element in the $i$the row and $j$th column of the matrix $I= \sigma^2I(\theta,G,{\cal P}),$ where the Fisher information matrix $I(\theta,G,{\cal P})$ is as defined in (\ref{efficientinfluencefunctiongeneral}) with $\theta =(\mu,\sigma)^T$ and ${\cal G} = \{g(\cdot)\}.$ Some computation shows $I_{11}= \int (g'/g)^2 g, I_{12}=I_{21}=\int x(g'/g(x))^2 g(x) dx,$ and $I_{22}= \int [xg'/g(x) + 1]^2 g(x) dx$ exist and are finite; cf. Section I.2.3 of \cite{Hajek67}.

We consider the submodel with the coefficient of variation known to be equal to a given constant $c=\sigma/\mu$ and with $\nu=\mu$ the parameter of interest.
Since in a parametric model the model itself is always least favorable, the conditions of Theorem \ref{onestep} are satisfied and the estimator $\hat{\nu}_n = \hat{\mu}_n$ of $\mu$ from (\ref{efficientestimatornu}) with $\bar{\nu}_n=\bar{\mu}_n,\,\hat{\theta}_n=(\bar{\mu}_n,\bar{\sigma}_n)^T,$ and ${\hat I}_n = {\bar \sigma}_n^{-2}I$ is efficient and some computation shows
\begin{equation}\label{cknowngeneral}
\hat{\mu}_n = \left(I_{11} + 2cI_{12} + c^2 I_{22}\right)^{-1}\left[\left(I_{11} + cI_{12}\right) \bar{\mu}_n + \left(I_{12} + cI_{22}\right) \bar{\sigma}_n \right].
\end{equation}
In case the density $g(\cdot)$ is symmetric around 0, the Fisher information matrix is diagonal and $\hat{\mu}_n$ from (\ref{cknowngeneral}) becomes
\begin{equation}\label{cknownsymmetric}
\hat{\mu}_n = \left(I_{11} + c^2 I_{22}\right)^{-1}\left[I_{11}\bar{\mu}_n + cI_{22} \bar{\sigma}_n \right].
\end{equation}
In the normal case with $g(\cdot)$ the standard normal density $\hat{\mu}_n$ reduces to
\begin{equation}\label{cknownnormal}
\hat{\mu}_n = (1 + c^2)^{-1}\left[\bar{\mu}_n + 2c\bar{\sigma}_n \right]
\end{equation}
with $\bar{\mu}_n$ and $\bar{\sigma}_n$ equal to e.g. the sample mean and the sample standard deviation, respectively; cf. \cite{Khan68}, \cite{GH76}, and \cite{Khan15}.
\end{exam}

\begin{exam}\label{Gaussiancopula}\textbf{Gaussian copula models}

Let $$\mathbf{X}_1=(X_{1,1},\ldots,X_{1,m})^T,\ldots,\mathbf{X}_n=(X_{n,1},\ldots,X_{n,m})^T$$ be i.i.d. copies of $\mathbf{X}=(X_1,\ldots,X_m)^T$. For $i=1,\ldots,m$, the marginal distribution function of $X_i$ is continuous and will be denoted by $F_i.$ It is assumed that $(\Phi^{-1}(F_1(X_1)), \dots, \Phi^{-1}(F_m(X_m)))^T$ has an $m$-dimensional normal distribution with mean 0 and positive definite correlation matrix $C(\theta),$ where $\Phi$ denotes the one-dimensional standard normal distribution function. Here the parameter of interest $\theta$ is the vector in $\mathbb{R}^{m(m-1)/2}$ that summarizes all correlation coefficients $\rho_{rs},\, 1\leq r < s \leq m$. We will set this general Gaussian copula model as our semiparametric starting model $\cal P$, i.e.,
\begin{equation}\label{generalGaussiancopula}
{\cal P}=\{P_{\theta,G}\ :\ \theta=(\rho_{12},\dots,\rho_{(m-1)m})^T\ , G =(F_1(\cdot),\dots,F_m(\cdot))\in {\cal G}\}.
\end{equation}
The unknown continuous marginal distributions are the nuisance parameters collected as $G \in {\cal G}.$

Theorem 3.1 of \cite{KW97} shows that the normal scores rank correlation coefficient is semiparametrically efficient in $\cal P$ for the 2-dimensional case with normal marginals with unknown variances constituting a least favorable parametric submodel. As \cite{Hoff14} explain at the end of their Section 1 and in their Section 4, their Theorem 4.1 proves that normal marginals with unknown, possibly unequal variances constitute a least favorable parametric submodel, also for the general $m$-dimensional case. Since the maximum likelihood estimators are efficient for the parameters of a multivariate normal distribution, the sample correlation coefficients are efficient for estimation of the correlation coefficients based on multivariate normal observations. But each sample correlation coefficient and hence its efficient influence function involve only two components of the multivariate normal observations. Apparently, the other components of the multivariate normal observations carry no information about the value of the respective correlation coefficient. Effectively, for each correlation coefficient we are in the 2-dimensional case and invoking again Theorem 3.1 of \cite{KW97} we see that also in the general $m$-dimensional case the normal scores rank correlation coefficients are semiparametrically efficient. They are defined as
\begin{equation}\label{normalscores}
\hat\rho_{rs}^{(n)}=\frac{\frac1n\sum\limits_{j=1}^{n}{\Phi^{-1}\left(\frac n{n+1}\mathbb{F}_r^{(n)}(X_{j,r})\right)\Phi^{-1}\left(\frac n{n+1}\mathbb{F}_s^{(n)}(X_{j,s})\right)}}{\frac1n\sum\limits_{j=1}^{n}{\left[\Phi^{-1}\left(\frac j{n+1}\right)\right]^2}}
\end{equation}
with $\mathbb{F}_r^{(n)}$ and $\mathbb{F}_s^{(n)}$ being the marginal empirical distributions of $F_r$ and $F_s$, respectively, $1\leq r < s \leq m.$
The Van der Waerden or normal scores rank correlation coefficient $\hat\rho_{rs}^{(n)}$ from (\ref{normalscores}) is a semiparametrically efficient estimator of $\rho_{rs}$ with efficient influence function
\begin{align}
\tilde\ell_{\rho_{rs}}(X_r,X_s)&=\Phi^{-1}\left(F_r(X_r)\right)\Phi^{-1}\left(F_s(X_s)\right)\\
&\ \ \ -\tfrac 12 \rho_{rs} \left\{\left[\Phi^{-1}\left(F_r(X_r)\right)\right]^2
+\left[\Phi^{-1}\left(F_s(X_s)\right)\right]^2\right\}\notag.
\end{align}
This means that
\begin{equation}\label{estimatorthetaexc}
\hat\theta_n=(\hat\rho_{12}^{(n)},\dots,\hat\rho_{(m-1)m}^{(n)})^T
\end{equation}
efficiently estimates $\theta$ with efficient influence function
\begin{equation}\label{estimatorinfthetaexc}
\tilde\ell(\mathbf{X}; \theta,G,{\cal P})=(\tilde\ell_{\rho_{12}}(X_1,X_2),\ldots,\tilde\ell_{\rho_{(m-1)m}}(X_{m-1},X_m))^T.
\end{equation}

\begin{subexam}\label{exchangeable} \textbf{Exchangeable Gaussian copula}

The exchangeable $m$-variate Gaussian copula model
\begin{equation}\label{exchangeablemodel}
{\cal Q} = \left\{ P_{{\bf 1}_k \rho,G}\,:\, \rho \in (-1/(m-1),1),\ G\in {\cal G}\right\} \subset {\cal P}
\end{equation}
is a submodel of the Gaussian copula model $\cal P$ with a one-dimensional parameter of interest $\nu=\rho.$ In this submodel all correlation coefficients have the same value $\rho.$ So, $\theta={\bf 1}_k \rho$ with ${\bf 1}_k$ indicating the vector of ones of dimension $k = m(m-1)/2.$
In order to construct an efficient estimator of $\rho$ within $\cal Q$ along the lines of Section \ref{EEPI}, in particular Remark \ref{remarklinear}, we first apply (\ref{nubarforlinear}) with $\alpha=0$ and $L= {\bf 1}_k$ to obtain the (natural) $\sqrt n$-consistent estimator
\begin{equation}\label{rhobar}
{\bar \rho}_n = {\bar \nu}_n = \frac 1k \sum_{r=1}^{m-1} \sum_{s=r+1}^m {\hat\rho}_{rs}^{(n)}.
\end{equation}
For $\theta={\bf 1}_k \rho$ we get by simple but tedious calculations (see the Supplementary Material)
\begin{equation}\label{covarianceinfluencefunctions}
E\tilde\ell_{\rho_{rs}}\tilde\ell_{\rho_{tu}}=
    \begin{cases}
      (1-\rho^2)^2 & \textrm{if}\quad |\{r,s\} \cap \{t,u\}|=2, \\
      \tfrac 12 (1-\rho)^2\rho(2+3\rho) & \textrm{if}\quad |\{r,s\} \cap \{t,u\}|=1, \\
      2(1-\rho)^2\rho^2 & \textrm{if}\quad |\{r,s\} \cap \{t,u\}|=0.\\
   \end{cases}
\end{equation}
It makes sense to estimate $I({\bf 1}_k,G,{\cal P})$ by substituting ${\bar \rho}_n$ for $\rho$ in (\ref{covarianceinfluencefunctions}), to compute the inverse of the resulting matrix, and to choose this matrix as the estimator ${\hat I}_n.$ To this end, we note that for every pair $\{r,s\},\, 1 \leq r \neq s \leq m$, there are $2(m-2)$ pairs of $\{t,u\}$'s having one element in common and there are $\frac12(m-2)(m-3)$ pairs of $\{t,u\}$'s having no elements in common. Hence, the sum of the components of each column vector of $I^{-1}({\bf 1}_k \rho,G, {\cal P})$ is $(1-\rho)^2(1+(m-1)\rho)^2.$ Each matrix with the components of each column vector adding to 1 has the property that the sum of all row vectors equals the vector with all components equal to 1, and hence the components of each column vector of its inverse also add up to 1. This implies
\begin{equation*}
{\bf 1}_k^T {\hat I}_n = \left(1-{\bar \rho}_n \right)^{-2} \left(1+(m-1){\bar \rho}_n \right)^{-2} {\bf 1}_k^T
\end{equation*}
and hence by (\ref{efficientnuforlinear})
\begin{equation}\label{efficientestimatorrhoexchangeable}
\hat\rho_n=\left({\bf 1}_k^T {\hat I}_n {\bf 1}_k \right)^{-1} {\bf 1}_k^T {\hat I}_n {\hat \theta}_n = \frac 1k {\bf 1}_k^T {\hat \theta}_n = {m \choose 2}^{-1}\sum_{r=1}^{m-1}\sum_{s=r+1}^m {\hat\rho_{rs}^{(n)}} = {\bar \rho}_n
\end{equation}
attains the asymptotic information bound (cf. (\ref{informationboundnu}))
\begin{equation}\label{informationboundexchangeable}
\left({\bf 1}_k^T I\left({\bf 1}_k\rho,G,{\cal P}\right) {\bf 1}_k \right)^{-1} = {m \choose 2}^{-1} (1-\rho)^2(1+(m-1)\rho)^2.
\end{equation}
\cite{Hoff14} proved the efficiency of the pseudo-likelihood estimator for $\rho$ in dimension $m=4$. \cite{Segers14} extended this result to general $m$ and presented the efficient lower bounds for $m=3$ and $m=4$ in their Example 5.3. However, their maximum pseudo-likelihood estimator is not as explicit as our (\ref{efficientestimatorrhoexchangeable}).
\end{subexam}

\begin{subexam}\label{circular} \textbf{Four-dimensional circular Gaussian copula}

A particular, one-dimensional parameter type of four-dimensional circular Gaussian copula model has been studied by \cite{Hoff14} and \cite{Segers14}. It is defined by its correlation matrix
\begin{equation}\label{cilcularmatrix}
 \begin{pmatrix}
  1 & \rho & \rho^2 & \rho \\
  \rho & 1 & \rho & \rho^2 \\
  \rho^2  & \rho  & 1 & \rho  \\
  \rho & \rho^2 & \rho & 1
 \end{pmatrix}.
\end{equation}
Our semiparametric starting model $\cal P$ is the same as in (\ref{generalGaussiancopula}) with $m=4,$ but with the components of $\theta$ rearranged as follows
\begin{equation*}
\theta=(\rho_{12}\ ,\ \rho_{14}\ ,\ \rho_{23}\ ,\ \rho_{34}\ ,\ \rho_{13}\ ,\ \rho_{24})^T.
\end{equation*}
Now, with $f(\rho)=(\rho\ ,\ \rho\ ,\ \rho\ ,\ \rho\ ,\ \rho^2\ ,\ \rho^2)^T$ the present circular Gaussian submodel $\cal Q$ may be written as
\begin{equation*}
{\cal Q}=\{P_{f(\rho),G}\ :\ \rho\in(-\tfrac 13,1)\,,\  G\in {\cal G}\}.
\end{equation*}
In order to construct an efficient estimator of $\rho$ within $\cal Q$ along the lines of Theorem \ref{onestep}, we propose as a $\sqrt n$-consistent estimator of $\rho$
\begin{eqnarray}\label{rootnconsistentestimatorrho}
\lefteqn{{\bar \rho}_n = \tfrac 23 {\bar \rho}_{n,1} + \tfrac 13\, {\rm sign}\left({\bar \rho}_{n,1}\right){\bar \rho}_{n,2}, \nonumber } \\
&& {\bar \rho}_{n,1} = \tfrac 14 \left({\hat\rho}_{12}^{(n)} + {\hat\rho}_{14}^{(n)} + {\hat\rho}_{23}^{(n)} + {\hat\rho}_{34}^{(n)}\right),\ {\bar \rho}_{n,2} = \tfrac 12 \left(\sqrt{{\hat\rho}_{13}^{(n)}} + \sqrt{{\hat\rho}_{24}^{(n)}} \right).
\end{eqnarray}
As in (\ref{covarianceinfluencefunctions}) we get by simple but tedious calculations (see the Supplementary Material)
\begin{equation}\label{inversematrix}
I^{-1}(f(\rho),G,{\cal P}) = \tfrac12 \left( 1-\rho^2\right)^2
\end{equation}
\begin{equation*}
\begin{pmatrix}
2 & \rho^2 & \rho^2 & 2\rho^2 & \rho\left( 2+\rho^2\right)
& \rho\left( 2+\rho^2\right) \\
\rho^2 & 2 & 2\rho^2 & \rho^2 & \rho\left( 2+\rho^2\right)
 & \rho\left( 2+\rho^2\right) \\
\rho^2 & 2\rho^2 & 2 & \rho^2 & \rho\left( 2+\rho^2\right)
 & \rho\left( 2+\rho^2\right) \\
2\rho^2 & \rho^2 & \rho^2 & 2 & \rho\left( 2+\rho^2\right) & \rho\left( 2+\rho^2\right) \\
\rho\left( 2+\rho^2\right) & \rho\left( 2+\rho^2\right) &
\rho\left( 2+\rho^2\right) & \rho\left( 2+\rho^2\right) &
2\left( 1+\rho^2\right) ^2 & 4\rho^2 \\
\rho\left( 2+\rho^2\right) & \rho\left( 2+\rho^2\right) &
\rho\left( 2+\rho^2\right) & \rho\left( 2+\rho^2\right) &
4\rho^2 & 2\left( 1+\rho^2\right) ^2
\end{pmatrix},
\end{equation*}
which has inverse
\begin{equation}\label{matrix}
I(f(\rho),G,{\cal P})= \tfrac 12 \left( 1-\rho^2\right)^{-4}
\end{equation}
\begin{equation*}
\begin{pmatrix}
\rho ^4+2 & 3\rho^2 & 3\rho^2 & \rho ^4+2\rho ^2 & -\left(\rho^3 + 2\rho \right) & -\left(\rho^3 + 2\rho \right) \\
3\rho^2 & \rho ^4+2 & \rho ^4+2\rho ^2 & 3\rho ^2 & -\left(\rho^3 + 2\rho \right) & -\left(\rho^3 + 2\rho \right) \\
3\rho^2 & \rho^4+2\rho ^2 & \rho ^4+2 & 3\rho ^2 & -\left(\rho^3 + 2\rho \right) & -\left(\rho^3 + 2\rho \right) \\
\rho ^4+2\rho ^2 & 3\rho ^2 & 3\rho^2 & \rho ^4+2 & -\left(\rho^3 + 2\rho \right) & -\left(\rho^3 + 2\rho \right) \\
 -\left(\rho^3 + 2\rho \right) &  -\left(\rho^3 + 2\rho \right) & -\left(\rho^3 + 2\rho \right) & -\left(\rho^3 + 2\rho \right) & 2\frac{\rho ^6+\rho ^4+1}{\rho^4 + 1} &  2\frac{\rho ^6+2\rho ^2}{\rho^4 +1} \\
 -\left(\rho^3 + 2\rho \right) & - \left(\rho^3 + 2\rho \right) & -\left(\rho^3 + 2\rho \right) & -\left(\rho^3 + 2\rho \right) & 2\frac{\rho ^6+2\rho ^2}{\rho^4 +1} & 2\frac{\rho ^6+\rho ^4+1}{\rho^4 + 1}
\end{pmatrix}.
\end{equation*}
Substituting ${\bar \rho}_n$ into (\ref{matrix}) we obtain a $\sqrt n$-consistent estimator of $I(f(\rho),G,{\cal P}).$
In view of ${\dot f}(\rho)= (1, 1, 1, 1, 2\rho, 2\rho)^T$ we have
\begin{equation*}
{\dot f}^T(\rho)I(f(\rho),G,{\cal P})= \left(1-\rho^2\right)^{-3}
\left(1+\rho^2, 1+\rho^2, 1+\rho^2, 1+\rho^2, -2\rho, -2\rho \right).
\end{equation*}
Consequently the asymptotic lower bound for estimation of $\rho$ within $\cal Q$ equals
\begin{equation}\label{boundcircular}
\left[{\dot f}(\rho)^T I({f(\rho)},G,{\cal P}) {\dot f}(\rho)\right]^{-1} = \tfrac 14 \left(1-\rho^2\right)^2.
\end{equation}
Substituting ${\bar \rho}_n$ for $\rho$ we obtain as the efficient estimator from Theorem \ref{onestep}
\begin{equation}
{\hat\rho}_n = {\bar \rho}_n + \frac {1 + {\bar \rho}_n^2}{1 - {\bar \rho}_n^2}\left({\bar \rho}_{n,1}- {\bar \rho}_n \right) - \frac {{\bar \rho}_n}{1 - {\bar \rho}_n^2} \left( \tfrac 12 \left(\hat\rho_{13}^{(n)}+ \hat\rho_{24}^{(n)}\right) - {\bar \rho}_n^2 \right).
\end{equation}

\cite{Hoff14} have shown that the pseudo-likelihood estimator is not efficient in this case. \cite{Segers14} have established the asymptotic lower bound (\ref{boundcircular}) and have constructed an alternative, efficient, one-step updating estimator suggesting the pseudo-maximum likelihood estimator as the preliminary estimator.
\end{subexam}
\end{exam}

\begin{exam}\label{regression}\textbf{Partial spline linear regression}

Here the observations are realizations of i.i.d. copies of the random vector $X=(Y,Z^T,U^T)^T$ with $Y, Z,$ and $U$ 1-dimensional, $k$-dimensional, and $p$-dimensional random vectors with the structure
\begin{equation}\label{partialspline}
Y = \theta^T Z + \psi(U) + \varepsilon,
\end{equation}
where the measurement error $\varepsilon$ is independent of $Z$ and $U,$ has mean 0, finite variance, and finite Fisher information for location, and where $\psi(\cdot)$ is a real valued function on ${\mathbb R}^p.$ \cite{Schick93} calls this partly linear additive regression, \cite{BKRW93} mention it as partial spline regression, whereas \cite{Cheng15} are talking about the partial smoothing spline model. Under the regularity conditions of his Theorem 8.1 \cite{Schick93} presents an efficient estimator of $\theta$ and a consistent estimator of $I(\theta, G, {\cal P}).$ Consequently our Theorem \ref{onestep} may be applied directly in order to obtain an efficient estimator of $\nu$ in appropriate submodels with $\theta=f(\nu)$ without our construction of an estimator of $I(\theta, G, {\cal P})$ via characteristic functions. Note that for submodels with $\theta$ restricted to a linear subspace, $\theta = L \nu$ say, our approach is not needed, since the reparametrization $Y = \nu^T L^T Z + \psi(U) + \varepsilon$ brings the estimation problem back to its original (\ref{partialspline}).
\end{exam}

\begin{exam}\label{commonmean}\textbf{Multivariate normal with common mean}

Let $\cal G$ be the collection of nonsingular $k\times k$-covariance matrices and let the parametric starting model be the collection of nondegenerate normal distributions with mean vector $\theta$ and covariance matrix $\Sigma,$
\begin{equation}\label{normalmodel}
{\cal P} = \left\{ P_{\theta,\Sigma} \,:\, \theta \in {\mathbb R}^k,\ \Sigma \in {\cal G} \right\}.
\end{equation}
Efficient estimators of $\theta$ and $\Sigma$ are the sample mean ${\bar X}_n = n^{-1}\sum_{i=1}^n X_i$ and the sample covariance matrix ${\hat \Sigma}_n = (n-1)^{-1} \sum_{i=1}^n (X_i - {\bar X}_n)(X_i - {\bar X}_n)^T,$ respectively. Note that ${\bar X}_n$ attains the finite sample Cram\'er-Rao bound and the asymptotic information bound with $I(\theta, \Sigma, {\cal P}) = \Sigma^{-1}.$

The parametric submodel we consider is
\begin{equation}\label{normalsubmodel}
{\cal Q} = \left\{ P_{{\bf 1}_k \mu,\Sigma} \,:\, \mu \in {\mathbb R},\ \Sigma \in {\cal G} \right\}.
\end{equation}
In view of (\ref{efficientnuforlinear}) and (\ref{boundcovariancematrix})
\begin{equation}\label{estimatorcommonmeans}
\hat\mu_n=\left(\mathbf{1}_k^T{\hat{\mathbf\Sigma}}_n^{-1} \mathbf{1}_k \right)^{-1}\mathbf{1}_k^T{\hat{\mathbf\Sigma}}_n^{-1}\mathbf{\bar X_n}
\end{equation}
is an efficient estimator of $\mu$ within $\cal Q$ that attains the asymptotic lower bound $\left(\mathbf{1}_k^T \Sigma^{-1} \mathbf{1}_k \right)^{-1}.$
In case the covariance matrix $\Sigma$ is diagonal with its variances denoted by $\sigma_1^2, \dots, \sigma_k^2,$ we are dealing with the Graybill-Deal model as presented by \cite{vanEeden06} on her page 88. With $\bar X_{i,n}=\frac1n\sum_{j=1}^n{X_{j,i}},\, S_{i,n}^2=\frac1n\sum_{j=1}^n (X_{j,i}-\bar X_{i,n})^2,$ and ${\hat \Sigma}_n = \textrm{diag}(S_{1,n}^2 ,\dots,S_{k,n}^2)$ we obtain the Graybill-Deal estimator
\begin{equation}\label{Graybill-Deal}
{\hat\mu}_n=\frac{\sum_{i=1}^k {\bar X_{i,n}/S_{i,n}^2}}{\sum_{i=1}^k{1/S_{i,n}^2}}
\end{equation}
with asymptotic lower bound $\left(\mathbf{1}_k^T \Sigma^{-1} \mathbf{1}_k \right)^{-1} = 1/\sum_{i=1}^k{1/\sigma_i^2}.$
\end{exam}

\begin{exam}\label{label}\textbf{Restricted maximum likelihood estimator}

Maximum likelihood estimation of the generalized linear model under linear restrictions on the parameters is done in \cite{Nyquist91} via an iterative procedure using a penalty function. \cite{Kim95} introduce the restricted EM algorithm for maximum likelihood estimation under linear restrictions. Our approach as described in Remark \ref{remarklinear} with ${\hat \theta}_n$ a(n unrestricted) maximum likelihood estimator avoids such iterative procedures.
\end{exam}

\appendix

\section{Additional Proofs}\label{TP}

In this appendix proofs will be presented of (\ref{efficientestimatorfnurestrictions}) and (\ref{informationboundrestriction}).

Since ${\hat I}_n$ has been chosen to be symmetric and positive definite, $x^T {\hat I}_n y,\, x,y \in {\mathbb R}^k,$ is an inner product on ${\mathbb R}^k.$ Define the $k \times k$-matrices $\Pi_{n,L}$ and $\Pi_{n,R}$ by
\begin{eqnarray}\label{projectionmatrices}
\Pi_{n,L} = L\left(L^T {\hat I}_n L \right)^{-1}L^T {\hat I}_n, \nonumber\\
\Pi_{n,R} = {\hat I}_n^{-1} R\left(R^T {\hat I}_n^{-1} R \right)^{-1}R^T.
\end{eqnarray}
With the above inner product these matrices are projection matrices on the linear subspaces spanned by the columns of $L$ and ${\hat I}_n^{-1}R,$ respectively. Indeed, $\Pi_{n,L}\Pi_{n,L} = \Pi_{n,L},\ \Pi_{n,R}\Pi_{n,R} = \Pi_{n,R},\ (x - \Pi_{n,L}x)^T {\hat I}_n \Pi_{n,L}x = 0,\, x \in {\mathbb R}^k,\ (y-\Pi_{n,R}y)^T {\hat I}_n \Pi_{n,R}y = 0,\, y \in {\mathbb R}^k,\ \Pi_{n,L}Lx = Lx,\, x \in {\mathbb R}^d,$ and $\Pi_{n,R}{\hat I}_n^{-1}Ry = {\hat I}_n^{-1}Ry,\, y \in {\mathbb R}^{k-d}$ hold.
The linear subspaces spanned by the columns of $L$ and ${\hat I}_n^{-1}R$ have dimensions $d$ and $k-d,$ respectively, since the matrices $(L,R)$ and ${\hat I}_n$ are nonsingular. Moreover, these linear subspaces are orthogonal in view of $L^T {\hat I}_n {\hat I}_n^{-1} R = L^T R =0.$ This implies
\begin{equation}\label{orthogonalprojections}
\Pi_{n,L}x + \Pi_{n,R}x = x,\quad x \in {\mathbb R}^k.
\end{equation}
Combining (\ref{projectionmatrices}), (\ref{orthogonalprojections}), and (\ref{efficientestimatorfnuforlinear}) we obtain (\ref{efficientestimatorfnurestrictions}) and, by the consistency of ${\hat I}_n,$ (\ref{informationboundrestriction}).

\section*{Acknowledgements}
We would like to thank Raymond Veldhuis for inspiring us to study problems with structured correlation matrices, which triggered the reported research, and Constance van Eeden for references.

\section*{Supplementary Material}
Computations needed for (\ref{covarianceinfluencefunctions}) and (\ref{inversematrix}) are collected as supplementary material.


\clearpage
\begin{frontmatter}
\begin{center}
\textbf{\LARGE Supplementary Material For "Semiparametrically Efficient Estimation of Constrained Euclidean Parameters"}
\end{center}
\end{frontmatter}
\setcounter{page}{1}
\makeatletter

In this supplement we present the computational details for \eqref{covarianceinfluencefunctions} and \eqref{inversematrix} presented in Example \ref{Gaussiancopula}. Since our computations will be based on fourth moments of multivariate normal random variables, we consider
\begin{equation*}
Z=
\begin{pmatrix}
Z_a\\Z_b\\Z_c\\Z_d
\end{pmatrix}
\sim
N\left(
\begin{pmatrix}
0\\0\\0\\0
\end{pmatrix}
,
\begin{pmatrix}
1&\rho_{ab}&\rho_{ac}&\rho_{ad}\\
\rho_{ba}&1&\rho_{bc}&\rho_{bd}\\
\rho_{ca}&\rho_{cb}&1&\rho_{cd}\\
\rho_{da}&\rho_{db}&\rho_{dc}&1\\
\end{pmatrix}
\right).
\end{equation*}
The following fourth moments of $Z$ can be obtained by straightforward computations:
\begin{itemize}
\item $E(Z_a^4)=3$
\item $E(Z_a^3Z_b)=3\rho_{ab}$
\item $E(Z_a^2Z_b^2)=1+2\rho_{ab}^2$
\item $E(Z_a^2Z_bZ_c)=\rho_{bc}+2\rho_{ab}\rho_{ac}$
\item $E(Z_aZ_bZ_cZ_d)=\rho_{ab}\rho_{cd}+\rho_{ac}\rho_{bd}+\rho_{ad}\rho_{bc}.$
\end{itemize}
For every $i,j=1,\dots,{n\choose2}$ let $M_{ij}$ be the element in the $i$-th row and $j$-th column of the efficient lower bound $I^{-1}(\theta,G,{\cal P}).$ Because of $\theta_i=\rho_{ab},\,\theta_j=\rho_{cd}$ for some $a,b,c,$ and $d,$ we have
\begin{equation*}
M_{ij}=E\left(Z_aZ_b-\tfrac12\rho_{ab}\left[Z_a^2+Z_b^2\right]\right)\left(Z_cZ_d-\tfrac12\rho_{cd}\left[Z_c^2+Z_d^2\right]\right).
\end{equation*}
We have three cases:
\begin{itemize}
\item $|\{a,b\} \cap \{c,d\}|=2$
\begin{eqnarray*}
\lefteqn{M_{ii} = E\left( Z_{a}Z_{b} - \tfrac 12 \rho_{ab} \left[Z_{a}^2+Z_{b}^2\right] \right)^2} \\
&& = E\left( Z_{a}^2Z_{b}^2\right) -\rho_{ab}E\left( Z_{a}^{3}Z_{b}+Z_{b}^{3}Z_{a}\right) +\tfrac 14 \rho_{ab}^2 E\left(Z_{a}^{4}+2Z_{a}^2Z_{b}^2+Z_{b}^{4}\right)  \\
&& = \left( 1+2\rho_{ab}^2\right) -\rho_{ab}\left( 3\rho_{ab}+3\rho_{ab}\right) +\tfrac 14 \rho_{ab}^2\left( 3+2\left[ 1+2\rho_{ab}^2
\right] +3\right)  \\
&& = \left( 1-\rho_{ab}^2\right)^2
\end{eqnarray*}

\item $|\{a,b\} \cap \{c,d\}|=1$ ({\it without lost of generality} assume $d=a$)
\begin{eqnarray*}
\lefteqn{M_{ij} = E\left( Z_{a}Z_{b}-\tfrac 12 \rho_{ab} \left[
Z_{a}^2+Z_{b}^2\right] \right) \left( Z_{a}Z_{c}-\tfrac 12 \rho_{ac}
\left[ Z_{a}^2+Z_{c}^2\right] \right) } \\
&& = E\left( Z_{a}^2Z_{b}Z_{c}\right) -\tfrac 12\rho_{ab} E\left(
Z_{a}^{3}Z_{c}+Z_{b}^2Z_{a}Z_{c}\right)  \\
&&\quad -\tfrac 12 \rho_{ac} E\left( Z_{a}^{3}Z_{b}+Z_{c}^2Z_{a}Z_{b}\right) \\
&&\quad +\tfrac 14 \rho_{ab}\rho_{ac} E\left(
Z_{a}^{4}+Z_{a}^2Z_{b}^2+Z_{a}^2Z_{c}^2+Z_{b}^2Z_{c}^2\right)  \\
&& = \left( \rho_{bc}+2\rho_{ab}\rho_{ac}\right) -\tfrac 12 \rho_{ab}
\left( 3\rho_{ac}+\left[ \rho_{ac}+2\rho_{ab}\rho_{bc}\right] \right) \\
&&\quad -\tfrac 12 \rho_{ac} \left( 3\rho_{ab}+\left[ \rho_{ab}+2\rho_{ac}\rho_{bc}\right] \right) \\
&&\quad +\tfrac 14 \rho_{ab}\rho_{ac} \left( 3+\left[ 1+2\rho_{ab}^2\right] +
\left[ 1+2\rho_{ac}^2\right] +\left[ 1+2\rho_{bc}^2\right] \right)  \\
&& = \tfrac 12 \left( 1- \rho_{ab}^2 - \rho_{ac}^2 \right) \left(2\rho_{bc} - \rho_{ab}\rho_{ac}\right) + \tfrac 12 \rho_{ab} \rho_{ac} \rho_{bc}^2
\end{eqnarray*}
\item $|\{a,b\} \cap \{c,d\}|=0$
\begin{eqnarray*}
\lefteqn{ M_{ij} = E\left( Z_{a}Z_{b}-\tfrac 12 \rho_{ab}\left[
Z_{a}^2+Z_{b}^2\right] \right) \left( Z_{c}Z_{d}-\tfrac 12 \rho_{cd}
\left[ Z_{c}^2+Z_{d}^2\right] \right) } \\
&& = E\left( Z_{a}Z_{b}Z_{c}Z_{d}\right) -\tfrac 12 \rho_{ab} E\left(
Z_{a}^2Z_{c}Z_{d}+Z_{b}^2Z_{c}Z_{d}\right)  \\
&&\quad -\tfrac 12\rho_{cd} E\left(
Z_{c}^2Z_{a}Z_{b}+Z_{d}^2Z_{a}Z_{b}\right)  \\
&&\quad +\tfrac 14 \rho_{ab}\rho_{cd} E\left(
Z_{a}^2Z_{c}^2+Z_{b}^2Z_{c}^2+Z_{a}^2Z_{d}^2+Z_{b}^2Z_{d}^2
\right)  \\
&& = \rho_{ab}\rho_{cd}+\rho_{ac}\rho_{bd}+\rho_{ad}\rho_{bc} -\tfrac 12 \rho_{ab} \left( \left[ \rho_{cd}+2\rho_{ac}\rho_{ad}\right] +\left[ \rho_{cd}+2\rho_{bc}\rho_{bd}\right] \right)  \\
&&\quad -\tfrac 12 \rho_{cd} \left( \left[ \rho_{ab}+2\rho_{ac}\rho_{bc}\right]
+\left[ \rho_{ab}+2\rho_{ad}\rho_{bd}\right] \right)  \\
&&\quad +\tfrac 14 \rho_{ab}\rho_{cd} \left( \left[ 1+2\rho_{ac}^2\right] +
\left[ 1+2\rho_{bc}^2\right] +\left[ 1+2\rho_{ad}^2\right] +\left[
1+2\rho_{bd}^2\right] \right)  \\
&& = \rho_{ac}\rho_{bd}+\rho_{ad}\rho_{bc} -\left( \rho_{ab}\rho_{ac}\rho_{ad}+\rho_{ba}\rho_{bc}\rho_{bd}+\rho_{ca}\rho_{cb}\rho_{cd}+\rho_{da}\rho_{db}\rho_{dc}\right)  \\
&&\quad +\tfrac 12 \rho_{ab}\rho_{cd} \left( \rho_{ac}^2+\rho_{bc}^2+\rho_{ad}^2+\rho_{bd}^2\right)
\end{eqnarray*}
\end{itemize}
Finally, substitution of the correlation structures in Subexample \ref{exchangeable} and Subexample \ref{circular} give \eqref{covarianceinfluencefunctions} and \eqref{inversematrix}, respectively.

\end{document}